\documentclass[11pt]{amsart}

\usepackage{amsmath,amssymb,latexsym,rawfonts}
\input{prepictex}
\input{pictex}
\input{postpictex}
\newcommand{\impl}{\longmapsto}
\newcommand{\LL}{{\it L}}

\newtheorem{theorem}{Theorem}[section]
\newtheorem{$^{*}$theorem}[theorem]{$^{*}$Theorem}
\newtheorem{lemma}[theorem]{Lemma}
\newtheorem{$^{*}$lemma}[theorem]{$^{*}$Lemma}
\newtheorem{definition}[theorem]{Definition}
\newtheorem{proposition}[theorem]{Proposition}
\newtheorem{$^{*}$proposition}[theorem]{$^{*}$Proposition}

\newtheorem{remark}[theorem]{Remark}

\numberwithin{equation}{theorem}
\begin{document}
\title{COMPACTNESS IN \ $L$-FUZZY TOPOLOGICAL SPACES}
\author{Joaqu\'\i n Luna-Torres $^{\lowercase{a}}$  \,\ and \,\ El\'\i as Salazar-Buelvas $^{\lowercase{b}}$}
\dedicatory{$^a$  Universidad Sergio arboleda\\
 $^b$  Universidad de Cartagena}
\email{$^a$jluna@ima.usergioarboleda.edu.co}
\email{$^b$elisalazar31@hotmail.com}
\subjclass[2000]{54A40, 54D30, 06B23}
\keywords{GL-monoid; Co-GL-monoid; L-fuzzy topological space; L-fuzzy filter; L-fuzzy neighborhood system; L-fuzzy interior operator; Adherent point; Tychonoff theorem.}
\footnote{Partially supported by Universidad de Cartagena}
\begin{abstract}
We give a definition of compactness in $L$-fuzzy topological spaces and provide a characterization of compact $L$-fuzzy topological spaces, where $L$ is a complete quasi-monoidal lattice with some additional structures, and we present a version of Tychonoff's theorem within the category of $L$-fuzzy topological spaces.
\end{abstract}
\maketitle
\baselineskip=1.7\baselineskip
\section*{0. Introduction}
Over the years, a number of descriptions of compactness of fuzzy topological spaces has appeared. As mathematicians sought to general topology in various ways using the concept of fuzzy subsets of an ordinary set, it is not surprising that searches for such properties were obtained with different degrees of success, depending on the structure of the
underlying lattice\ $L$.\ One of the main goals of any theory of compactness is to investigate the problem to which extent the formulation of a theorem of Tychonoff is possible. The aim of this paper is to present a version  of Tychonoff's theorem  within the category of $L$-fuzzy topological spaces, when the underlying lattice\ $L$\ is a cqm-lattice with some additional structures. Following P. T. Johnstone (\cite{Jhst}), within the text of the paper, those propositions, lemmas and theorems whose proofs require Zorn's lemma are distinguished by being marked with an asterisk.
\newline
The paper is organized as follows: After some lattice-theoretical prerequisites, where we briefly recall the concept of a cqm-lattice, in section  $2$ we shall present the concept of  compact $L$-fuzzy topological spaces. Finally, in section $3$ we present a proof of the Tychonoff's theorem.
\section{From Lattice Theoretic Foundations}
Let $(L, \leq)$ be a complete, infinitely distributive lattice, i.e. $(L, \leq)$ is a partially ordered set such that for every subset $A\subset L$ the join $\bigvee A$ and the meet $\bigwedge A$ are defined, moreover $(\bigvee A) \wedge \alpha = \bigvee \{ a\wedge \alpha) \mid a \in A \}$ and\linebreak $(\bigwedge A) \vee \alpha = \bigwedge \{a\vee \alpha) \mid a \in A \}$
for every $\alpha \in L$.  In particular, $\bigvee L =: \top$ and $\bigwedge L =: \bot$ are respectively the universal upper and the universal lower bounds in $L$.
We assume that $\bot \ne \top$, i.e.\ $L$ has at least two elements.

\subsection{$\mathbf cqm-$lattices} 
The definition of complete quasi-monoidal lattices introduced by  E. Rodabaugh in \cite{SER} is the following:

A $cqm-$lattice (short for complete quasi-monoidal lattice) is a
triple\linebreak $(l,\leqslant,\otimes)$ provided with the
following properties
\begin{enumerate}
\item[(1)] $(L,\leqslant)$ is a complete lattice with upper bound $\top$ and lower bound $\bot$.
\item[(2)] $\otimes: L\times L : \rightarrow L$ is a binary operation satisfying the following axioms:
\begin{enumerate}
\item $\otimes$ is isotone in both arguments, i.e. $\alpha_1 \leqslant\alpha_2,\,\ \beta_1\leqslant \beta_2$ implies
$\alpha_1\otimes\beta_1\leqslant \alpha_2\otimes\beta_2$; \item  $\top$ is idempotent, i.e. $\top \otimes\top = \top$.
\end{enumerate}
\end{enumerate}

\subsection{$GL-$monoids}
A $GL-$monoid (see \cite{Ho91}, \cite{Ho92}, \cite{Ho94}) is a complete
lattice enriched with a further binary operation $\otimes$, i.e.\ a triple $(L, \leq, \otimes)$ such that:
\begin{enumerate}
\item[(1)]
$\otimes$ is isotone, i.e.\ $\alpha \leq \beta$ implies $\alpha \otimes \gamma \leq \beta \otimes \gamma$, $\forall \alpha, \beta, \gamma \in \LL$;
\item[(2)]
$\otimes$ is commutative, i.e.\ $\alpha \otimes \beta = \beta \otimes \alpha$, $\forall \alpha, \beta \in \LL$;
\item[(3)]
$\otimes$ is associative, i.e.\ $\alpha \otimes (\beta \otimes \gamma) = (\alpha \otimes \beta) \otimes \gamma$, $\forall \alpha, \beta, \gamma \in L$;
\item[(4)]
$(L,\leq,\otimes)$ is integral, i.e.\ $\top$ acts as the unity: $\alpha \otimes \top = \alpha$, $\forall \alpha \in \LL$;
\item[(5)]
$\bot$ acts as the zero element in $(L, \leq, \otimes)$, i.e.\ $\alpha\otimes \bot = \bot$, $\forall \alpha \in \LL$;
\item[(6)]
$\otimes$ is distributive over arbitrary joins, i.e.\ $\alpha \otimes (\bigvee_{\lambda} \beta_{\lambda}) = \bigvee_{\lambda} (\alpha \otimes \beta_{\lambda})$,
$\forall \alpha \in \LL, \forall \{ \beta_{\lambda} : \lambda \in I\} \subset \LL$;
\item[(7)]
$(L, \leq, \otimes)$ is divisible, i.e.\ $\alpha \leq \beta$ implies the existence of $\gamma \in L$ such that $\alpha = \beta \otimes \gamma$.
\end{enumerate}
It is well known that every $GL-$monoid is residuated, i.e.\ there exists a further binary operation ``$\impl$'' (implication) on $L$ satisfying the following condition:
$$\alpha \otimes \beta \leq \gamma \Longleftrightarrow \alpha \leq (\beta \impl \gamma) \qquad \forall \alpha, \beta, \gamma \in L.$$
Explicitly implication is given by
\[
\alpha \impl \beta = \bigvee \{ \lambda \in L \mid \alpha \otimes \lambda \leq \beta \}.
\]
Important examples of $GL$-monoids are Heyting algebras and $MV$-alg\-ebras. Namely, a {\it Heyting algebra} is $GL$-monoid of the type $(L,\leq,\wedge,\vee,\wedge)$ (i.e.\ in case of a Heyting algebra $\wedge = \otimes$), cf.\ e.g.\ \cite{Jhst}.
A $GL$-monoid is called an {\it $MV$-algebra} if $(\alpha \impl \bot) \impl \bot = \alpha \quad \forall \alpha \in L$,
\cite{Ch58}, \cite{Ch59}, see also \cite[Lemma 2.14]{Ho94}. Thus in an $MV$-algebra an order reversing involution $^c: L \to L$ can be naturally defined by setting $\alpha^c := \alpha \impl \bot \quad \forall \alpha \in L$.
\newline
If $X$ is a set and $L$ is a $GL$-monoid, then the fuzzy powerset $L^X$ in an obvious way can be pointwise endowed with a structure of a $GL$-monoid. In particular the $L$-sets $1_X$ and $0_X$ defined by $1_X (x):= \top$ and $0_X (x) := \bot$ $\forall x \in X$ are respectively the universal  upper and lower bounds in $L^X$.
\subsection{Co-$GL$-monoids}
A co-$GL-$monoid is a complete lattice with a further binary operation $\oplus$, i.e.\ a triple $(L, \leqslant, \oplus)$ such that:
\begin{enumerate}
\item[(1)]
$\oplus $ is isotone, i.e.\ $\alpha \leqslant \beta$ implies $\alpha \oplus \gamma \leqslant \beta \oplus \gamma$, $\forall \alpha, \beta, \gamma \in \LL$;
\item[(2)]
$\oplus $ is commutative, i.e.\ $\alpha \oplus \beta = \beta \oplus \alpha$, $\forall \alpha, \beta \in \LL$;
\item[(3)]
$\oplus $ is associative, i.e.\ $\alpha \oplus (\beta \oplus \gamma) = (\alpha \oplus \beta) \oplus \gamma$, $\forall \alpha, \beta, \gamma \in L$;
\item[(4)]
$(L,\leq,\oplus)$ is co-integral, i.e.\ $\bot$ acts as the unity: $\alpha \oplus \bot = \alpha$, $\forall \alpha \in \LL$;
\item[(5)]
$\top$ acts as the co-zero element in $(L, \leqslant, \oplus)$, i.e.\ $\alpha\oplus \top = \top$, $\forall \alpha \in \LL$;
\item[(6)]
$\oplus$ is distributive over arbitrary meets, i.e.\ $\alpha \oplus (\bigwedge_{\lambda} \beta_{\lambda}) = \bigwedge_{\lambda} (\alpha \oplus \beta_{\lambda})$, $\forall \alpha \in \LL, \forall \{ \beta_{\lambda} : \lambda \in I \} \subset \LL$;
\item[(7)]
$(L, \leqslant, \oplus)$ is co-divisible, i.e.\ $\alpha \leqslant \beta$ implies the existence of $\gamma \in L$ such that $\alpha \oplus \gamma= \beta $.
\end{enumerate}
Every co-$GL-$monoid is co-residuated, i.e.\ there exists a further binary operation ``$\leftarrowtail$'' (co-implication) on $L$ satisfying the following condition:
\[
\alpha \leftarrowtail\beta \leqslant \gamma \Longleftrightarrow \alpha \leqslant (\beta \oplus \gamma) \qquad \forall \alpha, \beta, \gamma \in L.
\]
Explicitly, coimplication is given by
\[
\alpha \leftarrowtail \beta = \bigwedge \{\lambda \in L \mid \alpha  \leqslant \beta\oplus\lambda \}.
\]
If $X$ is a set and $L$ is a co-$GL$-monoid, then the fuzzy powerset $L^X$ in an obvious way can be pointwise endowed with a structure of a co-$GL$-monoid.
\begin{remark} In this paper we will use the particular case in which  $\oplus=\lor$ and the co-implication is
\[
\alpha \vartriangleright\beta =\bigwedge \{\lambda \in L \mid \alpha  \leqslant \beta\lor\lambda \}.
\]
\end{remark}
In the sequel $(L,\leqslant, \otimes)$ denotes an arbitrary $GL$-monoid and  $(L,\leqslant, \lor)$ denotes an co-$GL$-monoid.
\section{Compactness and $L$-fuzzy Topological Spaces}
In this section we recall some topological concepts (as in \cite{HS}) from the theory of $L$-fuzzy topological spaces.
\begin{definition}
Given\,\ $(L,\leqslant, \otimes)$\,\ a cqm-lattice and\,\ $X$\,\ a non-empty set; an\,\ {\it $L$-fuzzy topology} on\,\ $X$\,\ is a map\,\ $\tau: L^X \to L$\,\ satisfying the following axioms:
\begin{enumerate}
\item[o1.] $\tau(1_X)=\top$,
\item[o2.] For \,\ $f,g \in L^X$,\,\ $\tau(f)\otimes \tau(g)\leqslant \tau(f\otimes g)$,
\item[o3.] For every subset\,\ $\{ f_{\lambda}\}_{\lambda \in I}$\,\ of\,\ $L^X$\,\ the inequality
\[
\bigwedge_{\lambda \in I}\tau(f_{\lambda})\leqslant \tau(\bigvee_{\lambda \in I}f_{\lambda})
\]
holds.
\end{enumerate}
If\,\ $\tau$\,\ is an\,\ $L$-fuzzy topology on\,\ $X$,\,\ the pair\,\ $(X,\tau)$\,\ is an\,\ $L$-fuzzy topological space.
\end{definition}
If the set\,\ $X$\,\ is non-empty, then \,\ $L^X$\,\ consists at least of two elements. In particular, the universal lower bound in\,\ $(L^X,\leqslant)$\,\ is given by\,\ $0_X$.\,\ When we take the empty subset of\,\ $L^X$\,\ and we apply the axiom o3, we obtain
\begin{enumerate}
\item[o1'.] $\tau(0_X)= \top $.
\end{enumerate}
\begin{definition}
Given\,\ $(X,\tau)$\,\ and\,\ $(Y,\eta)$\,\ $L$-fuzzy topological spaces, a map\,\ $\phi: X\to Y$\,\ is\,\ {\it $L$-fuzzy continuous} iff for all \ $g\in L^Y$,\,\ $\phi$\,\ satisfies
\[
\eta(g)\leqslant \tau(g\circ \phi).
\]
\end{definition}
In this way, we have the category  $\mathbf {L\mbox{-}FTop}$ where the objects are the\linebreak $L$-fuzzy topological spaces and the morphisms are the \,\ $L$-fuzzy continuous maps.
\newline
On the set $\mathfrak T_L(X)$ of all L-fuzzy topologies on $X$ we have a partial ordering:
\[
\tau_1\leqslant \tau_2 \,\,\,\ \Leftrightarrow \,\,\,\ \tau_1(f)\leqslant \tau_2(f),\,\,\,\,\,\ \forall f\in L^X.
\]
\begin{remark}
On the set $L^X\times L$ we introduce a partial ordering by
\[
(f,\alpha)\preccurlyeq (g,\beta) \,\,\ \Leftrightarrow \,\,\ \left(f\leqslant g\,\ \text{and}\,\  \beta\leqslant \alpha\right)
\]
\end{remark}
\begin{lemma}
Let $L$ be a complete lattice such that $(L, \leqslant, \otimes)$ is a $GL-$monoid and $(L, \leqslant, \lor)$ is a co-$GL-$monoid, then  $(L^X\times L, \preccurlyeq, \boxtimes )$, where
\[
(f,\alpha)\boxtimes (g,\beta):= (f\otimes g, \alpha\lor \beta),\,\,\,\,\,\,\,\ \forall (f,\alpha),(g,\beta)\in L^X\times L,
\]
 is a $GL-$monoid.
\end{lemma}
\begin{proof}
It is straightforward to check that for $\{(f_{\lambda},\alpha_{\lambda})\}_{\lambda\in I}\subseteq L^X\times L $:
\begin{enumerate}
\item[(i)] $\bigwedge_{\lambda\in I}(f_{\lambda},\alpha_{\lambda})=\left(\bigwedge_{\mbox{\tiny $\lambda\in I$}}f_{\lambda},\bigvee_{\mbox{\tiny $\lambda\in I$}}\alpha_{\lambda}\right)$.
\item[(ii)] $\bigvee_{\lambda\in I}(f_{\lambda},\alpha_{\lambda})=\left(\bigvee_{\mbox{\tiny $\lambda\in I$}}f_{\lambda},\bigwedge_{\mbox{\tiny $\lambda\in I$}}\alpha_{\lambda}\right)$.
\item[(iii)] $\boldsymbol{\top}=\left( 1_X,\bot\right)$ is the universal upper bound of $(L, \leqslant, \otimes)$.
\item[(iv)] $\boldsymbol{\bot}=\left( 0_X,\top\right)$ is the universal lower bound of $(L, \leqslant, \otimes)$.
\end{enumerate}
Consequently,  $\boxtimes$ is isotone, commutative and associative; $(L, \leqslant, \otimes)$ is integral with respect to $\boldsymbol{\top}$ and $\boldsymbol{\bot}$ acts as the zero element in $(L, \leqslant, \otimes)$.
\begin{enumerate}
\item[(v)] $\boxtimes$ is distributive over arbitrary joins, i.e. for all $ (f,\alpha)\in L^X\times L$,\ for all $\{(g_j,\beta_j)\mid j\in I\}\subseteq L^X\times L$,
\begin{align*}
(f,\alpha)\boxtimes\left[\bigvee_{j\in I}(g_j,\beta_j)\right]&= (f,\alpha)\boxtimes\left(\bigvee_{j\in I}g_j, \bigwedge_{j\in I}\beta_j \right)\\
&=\left( f\otimes\left(\bigvee_{j\in I}g_j\right),\alpha\lor\left(\bigwedge_{j\in I}\beta_j\right) \right)\\
&=\left( \bigvee_{j\in I}(f\otimes g_j),\bigwedge_{j\in I}(\alpha\lor\beta_j) \right)\\
&= \bigvee_{j\in I}\left(f\otimes g_j,\alpha\lor\beta_j \right)\\
&= \bigvee_{j\in I}\left[(f,\alpha)\boxtimes(g_j,\beta_j) \right].\\
\end{align*}
\item[(vi)] $(L^X\times L, \preccurlyeq, \boxtimes )$ is divisible, i.e. $(f,\alpha)\preccurlyeq (g,\beta)\Leftrightarrow f\leqslant g\,\ \text{and}\,\  \beta\leqslant \alpha $ implies the existence of $h\in L^X$ and $\gamma \in L$ such that $f=g\otimes h$ and $\beta\lor\gamma=\alpha$, in other words $(f,\alpha)=(g\otimes h, \beta \lor \gamma)=(g,\beta)\boxtimes (h,\gamma).$\qed
\end{enumerate}
\renewcommand{\qed}{}
\end{proof}
\begin{remark}
The residuation in $(L^X\times L, \preccurlyeq, \boxtimes )$ is the binary operation $\boldsymbol{\impl}$ given by
\begin{align*}
(f,\alpha)\boldsymbol{\impl} (g,\beta)&= \bigvee_{h\in L^X}\bigwedge_{\eta \in L}\left\{(h,\eta)\mid h\otimes f\leqslant g,\,\  \beta \leqslant \eta\lor\alpha \right\}\\
&=(f\impl g, \beta\rhd \alpha),\\
\end{align*}
which satisfies the following condition
\[
(f,\alpha)\boxtimes (g,\beta)\preccurlyeq (h,\gamma) \Leftrightarrow (f,\alpha)\preccurlyeq \left[(g,\beta)\boldsymbol{\impl} (h,\gamma)\right],
\]
$\forall (f,\alpha),(g,\beta),(h,\gamma)\in L^X\times L.$
It is also  straightforward to check that
\begin{itemize}
\item $(f,\alpha)\boxtimes\left[ (g,\beta)\boldsymbol{\impl} (h,\gamma)\right]\preccurlyeq\left[ (g,\beta)\boldsymbol{\impl}(f,\alpha)\boxtimes(h,\gamma)\right].$
\item $ \left[(g,\beta)\boldsymbol{\impl} (f,\alpha)\right] \boxtimes \left[(g,\beta)\boldsymbol{\impl} (h,\gamma)\right]\preccurlyeq\left[ (g,\beta)\boldsymbol{\impl}(f,\alpha)\boxtimes(h,\gamma)\right],$ whenever $\otimes$ will be idempotent.
\end{itemize}
\end{remark}
\begin{definition}[\bf{$L$-fuzzy filter}] 
Let $X$ be a set.  A map $\mathcal F:  L^X\times L\rightarrow L$ is called an $L$-fuzzy filter on $X$  if and only if $\mathcal F$ satisfies the following axioms:
\begin{enumerate}
\item[(FF0)] $\mathcal F(1_X, \alpha)= \top$.
\item[(FF1)] $(f,\alpha)\preccurlyeq (g,\beta)\Rightarrow  \mathcal F(f, \alpha)\leqslant \mathcal F(g, \beta)$.
\item[(FF2)] $\mathcal F(f, \alpha)\otimes \mathcal F(g, \beta)\leqslant \mathcal F(f\otimes g, \alpha\lor\beta)$.
\item[(FF3)] $\mathcal F(0_X, \alpha)= \bot$.
\end{enumerate}
\end{definition}
Let $\mathfrak F_{LF}(X)$ be the set of all $L$-fuzzy filters on $X$. On  $\mathfrak F_{LF}(X)$ we introduce a partial ordering \ $\curlyeqprec$\ by
\[
\mathcal F_1\curlyeqprec \mathcal F_2 \Leftrightarrow \mathcal F_1(f,\alpha)\leqslant \mathcal F_2(f,\alpha),\,\,\,\,\ \forall (f,\alpha)\in L^X\times L
\]
\begin{$^{*}$proposition}
The partially ordered set $(\mathfrak F_{LF}(X),\curlyeqprec)$ has maximal elements.
\end{$^{*}$proposition}
\begin{proof}
Referring to Zorn's lemma, it is sufficient to show  that every chain $\mathcal C$ in $\mathfrak F_{LF}(X)$ has an upper bound in $\mathfrak F_{LF}(X)$. For this purpose let us consider a non-empty chain $\mathcal C=\{\mathcal F_{\lambda}\mid \lambda \in I\} $. We define a map $\mathcal F_{\infty}: L^X\times L\rightarrow L$ by
\[
\mathcal F_{\infty}(f, \alpha)=\bigvee_{\lambda\in I}\mathcal F_{\lambda}(f, \alpha),
\]
and we show that $\mathcal F_{\infty}$ is an $L$-fuzzy filter on $X$. In fact
\begin{enumerate}
\item[(FF0)] $\mathcal F_{\infty}(1_X, \alpha)=\bigvee_{\lambda\in I}\mathcal F_{\lambda}(1_X, \alpha)=\bigvee_{\lambda\in I}\top=\top$.
\item[(FF1)] $(f,\alpha)\preccurlyeq (g,\beta)\Rightarrow  \mathcal F_{\infty}(f, \alpha)=\bigvee_{\lambda\in I}\mathcal F_{\lambda}(f, \alpha)\leqslant\bigvee_{\lambda\in I}\mathcal F_{\lambda}(g, \beta) =\mathcal F_{\infty}(g, \beta)$.
\item[(FF2)]
\begin{align*}
\mathcal F_{\infty}(f, \alpha)\otimes\mathcal F_{\infty}(g, \beta)&= \left(\bigvee_{\lambda\in I}\mathcal F_{\lambda}(f, \alpha)\right)\otimes \left(\bigvee_{\lambda\in I}\mathcal F_{\lambda}(g, \beta)\right)\\
&=\bigvee_{\lambda\in I}\left[\mathcal F_{\lambda}(f, \alpha)\otimes\mathcal F_{\lambda}(g, \beta)\right]\\
&\leqslant \bigvee_{\lambda\in I}\left[\mathcal F_{\lambda}(f\otimes g, \alpha\lor\beta)\right]\\
&=\mathcal F_{\infty}(f\otimes g, \alpha\lor\beta).
\end{align*}
\item[(FF3)]  $\mathcal F_{\infty}(0_X, \alpha)=\bigvee_{\lambda\in I}\mathcal F_{\lambda}(0_X, \alpha)=\bigvee_{\lambda\in I}\bot=\bot$.
\end{enumerate}
\end{proof}
\begin{definition}
A maximal element in $(\mathfrak F_{LF}(X),\curlyeqprec)$ is also called an $L$-fuzzy ultrafilter.
\end{definition}
\begin{$^{*}$proposition}\label{UF}
For every  $L$-fuzzy filter $\mathcal U:  L^X\times L\rightarrow L$  on $X$ the following assertions are equivalent
\begin{enumerate}
\item[(i)] $\mathcal U$ is an $L$-fuzzy ultrafilter.
\item[(ii)] $\mathcal U(f, \alpha)=\left(\mathcal U\left[(f,\alpha)\boldsymbol{\impl}(0_X,\rho)\right]\right)\impl\bot,\,\ \text{for all}\,\ (f,\alpha)\in L^X\times L,\,\\ \text {for all}\,\ \rho\leqslant\alpha\,\ \text{in} \,\ L. $
\end{enumerate}
\end{$^{*}$proposition}
\begin{proof} $(i)\Rightarrow(ii)$\newline
Because of $(FF2)$ and $(FF3)$  every $L$-fuzzy filter satisfies the condition
\begin{enumerate}
\item[(FF3')]  $\mathcal U(f, \alpha)\leqslant\left(\mathcal U\left[(f, \alpha)\boldsymbol{\impl} (0_X,\rho)\right]\right)\impl\bot,\,\ \text{for all}\,\ (f,\alpha)\in L^X\times L, \\ \text{for all}\,\ \rho\leqslant\alpha\,\ \text{in} \,\ L.$
\end{enumerate}
In order to verify $(i)\Rightarrow (ii)$ it is sufficient to show that the maximality of $\mathcal U$ implies
\[
\left(\mathcal U\left[(f, \alpha)\boldsymbol{\impl} (0_X,\rho)\right]\right)\impl\bot\leqslant \mathcal U(f, \alpha),\,\ \forall (f,\alpha)\in L^X\times L,\ \forall\rho\leqslant\alpha\ \text{in}\ L.
\]
For this purpose, we fix an element $(g,\beta)\in L^X\times L$, for that element we let $\mathcal G_{\beta}:= \left(\mathcal U\left[(g, \beta)\boldsymbol{\impl} (0_X,\rho)\right]\right)\impl\bot$ and define a map\ $\hat{\mathcal U}:  L^X\times L\rightarrow L$ \ by
\[
\hat{\mathcal U}(f,\alpha)=\mathcal U(f, \alpha)\bigvee\Bigl\{\mathcal U\left[(g, \beta)\boldsymbol{\impl} (f,\alpha)\right]\otimes \mathcal G_{\beta}\Bigr\}.
\]
We must show that $\hat{\mathcal U}$ is an $L$-fuzzy ultrafilter.
Firstly $\hat{\mathcal U}$ is an $L$-fuzzy filter: obviously
$\hat{\mathcal U}$ satisfies $(FF0)$. For the axiom $(FF1)$, from
the definition
\[
\hat{\mathcal U}(f,\alpha)=\mathcal U(f, \alpha)\bigvee\Bigl\{\mathcal U\left[(g, \beta)\boldsymbol{\impl} (f,\alpha)\right]\otimes \mathcal G_{\beta}\Bigr\}
\]
and
\[
\hat{\mathcal U}(h,\gamma)=\mathcal U(h, \gamma)\bigvee\Bigl\{\mathcal U\left[(g, \beta)\boldsymbol{\impl} (h,\gamma)\right]\otimes \mathcal G_{\beta}\Bigr\}.
\]
Now, for $(f,\alpha)\curlyeqprec (h,\gamma)$ we have that $\mathcal U(f,\alpha)\leqslant\mathcal U(h,\gamma)$,
moreover,
\begin{align*}
(g, \beta)\boldsymbol{\impl} (f,\alpha)&=\bigvee\bigl\{(k,\delta)\in L^X\times L\mid (g,\beta)\boxtimes (k, \delta)\curlyeqprec (f,\alpha)\bigr\}\\
&\leqslant \bigvee\bigl\{(k,\delta)\in L^X\times L\mid (g,\beta)\boxtimes (k, \delta)\curlyeqprec (h,\gamma)\bigr\}\\
&=(g, \beta)\boldsymbol{\impl} (h,\gamma),
\end{align*}
which implies that
\begin{align*}
\hat{\mathcal U}(f,\alpha)&=\mathcal U(f, \alpha)\bigvee\Bigl\{\mathcal U\left[(g, \beta)\boldsymbol{\impl} (f,\alpha)\right]\otimes \mathcal G_{\beta}\Bigr\}\\
&\leqslant\mathcal U(h, \gamma)\bigvee\Bigl\{\mathcal U\left[(g, \beta)\boldsymbol{\impl} (h,\gamma)\right]\otimes \mathcal G_{\beta}\Bigr\}\\
&=\hat{\mathcal U}(h,\gamma).
\end{align*}
For the axiom $(FF2)$, we must verify that
\[
\hat{\mathcal U}(f,\alpha)\otimes \hat{\mathcal U}(h,\gamma)\leqslant \hat{\mathcal U}(f\otimes h,\alpha\lor \gamma),
\]
In fact,\noindent
\begin{align*}
&\hat{\mathcal U}(f,\alpha)\otimes \hat{\mathcal U}(h,\gamma)\\
&=\bigl(\mathcal U(f, \alpha)\bigvee\bigl\{\mathcal U\left[(g, \beta)\boldsymbol{\impl} (f,\alpha)\right]\otimes \mathcal G_{\beta}\bigr\}\bigr)\\
&\otimes \bigl(\mathcal U(h, \gamma)\bigvee\bigl\{\mathcal U\left[(g, \beta)\boldsymbol{\impl} (h,\gamma)\right]\otimes \mathcal G_{\beta}\bigr\}\bigr)\\
&=\mathcal U(f, \alpha)\otimes\left[\mathcal U(h, \gamma)\bigvee\bigl\{\mathcal U\left[(g, \beta)\boldsymbol{\impl} (h,\gamma)\right]\otimes \mathcal G_{\beta}\bigr\}\right]\\
&\bigvee\left[\bigl\{\mathcal U\left[(g, \beta)\boldsymbol{\impl} (f,\alpha)\right]\otimes \mathcal G_{\beta}\bigr\}
\otimes\bigl\{\mathcal U(h, \gamma)\bigvee\bigl\{\mathcal U\left[(g, \beta)\boldsymbol{\impl} (h,\gamma)\right]\otimes \mathcal G_{\beta}\bigr\}\bigr\}\right]\\
&=\mathcal U(f, \alpha)\otimes \mathcal U(h, \gamma)\bigvee\left[\mathcal U(f, \alpha)\otimes\bigl\{\mathcal U\left[(g, \beta)\boldsymbol{\impl} (h,\gamma)\right]\otimes \mathcal G_{\beta}\bigr\}\right]\\
&\bigvee\left(\bigl\{\mathcal U\left[(g, \beta)\boldsymbol{\impl} (f,\alpha)\right]\otimes \mathcal G_{\beta} \bigr\}\otimes \mathcal U(h, \gamma)\right)\\
&\bigvee \left(\bigl\{ \mathcal U\left[(g, \beta)\boldsymbol{\impl} (f,\alpha)\right]\otimes \mathcal G_{\beta} \bigr\}
\otimes \bigl\{\mathcal U\left[(g, \beta)\boldsymbol{\impl} (h,\gamma)\right]\otimes \mathcal G_{\beta}\bigr\}\right)\\
&=\mathcal U(f, \alpha)\otimes \mathcal U(h, \gamma)\bigvee\left[\bigl\{\mathcal U(f, \alpha)\otimes\mathcal U\left[(g, \beta)\boldsymbol{\impl} (h,\gamma)\right]\bigr\}\otimes \mathcal G_{\beta}\right]\\
&\bigvee\left(\bigl\{\mathcal U(h, \gamma)\otimes\mathcal U\left[(g, \beta)\boldsymbol{\impl} (f,\alpha)\right]\bigr\} \otimes \mathcal G_{\beta} \right)\\
&\bigvee \left(\bigl\{ \mathcal U\left[(g, \beta)\boldsymbol{\impl} (f,\alpha)\right] \otimes \mathcal U\left[(g, \beta)\boldsymbol{\impl} (h,\gamma)\right]\bigr\}\otimes \mathcal G_{\beta}\right)\\
&\leqslant \mathcal U(f\otimes h, \alpha\lor\gamma)\bigvee\left[\mathcal U\left((f, \alpha)\boxtimes\left[(g, \beta)\boldsymbol{\impl} (h,\gamma)\right] \right)\otimes\mathcal G_{\beta}\right]\\
&\bigvee\left(\mathcal U\left((h,\gamma)\boxtimes\left[(g, \beta)\boldsymbol{\impl} (f,\alpha)\right]\right)\otimes\mathcal G_{\beta}\right)\\
&\bigvee\left(\mathcal U\left(\left[(g, \beta)\boldsymbol{\impl} (f,\alpha)\right]\boxtimes \left[(g, \beta)\boldsymbol{\impl} (h,\gamma)\right]\right)\otimes\mathcal G_{\beta}\right)\\
&\leqslant\mathcal U(f\otimes h, \alpha\lor\gamma)\bigvee\left[\mathcal U\left[(g,\beta)\boldsymbol{\impl}(f,\alpha)\boxtimes (h,\gamma)\right]\otimes\mathcal G_{\beta}\right]\\
&\bigvee\left[\mathcal U\left[(g,\beta)\boldsymbol{\impl}(f,\alpha)\boxtimes (h,\gamma)\right]\otimes\mathcal G_{\beta}\right]\\
&\bigvee\left[\mathcal U\left[(g,\beta)\boldsymbol{\impl}(f,\alpha)\boxtimes (h,\gamma)\right]\otimes\mathcal G_{\beta}\right]\\
&=\mathcal U(f\otimes h, \alpha\lor\gamma)\bigvee\left[\mathcal U\left[(g,\beta)\boldsymbol{\impl}(f,\alpha)\boxtimes (h,\gamma)\right]\otimes\mathcal G_{\beta}\right]\\
&=\hat{\mathcal U}(f\otimes h, \alpha\lor\gamma).
\end{align*}
\vspace{0.5cm}
In order to verify $(FF3)$, we have that
\begin{align*}
\hat{\mathcal U}(0_X,\alpha)&=\mathcal U(0_X, \alpha)\bigvee\Bigl\{\mathcal U\left[(g, \beta)\boldsymbol{\impl} (0_X,\alpha)\right]\otimes \mathcal G_{\beta}\Bigr\}\\
&=\bot\lor\Bigl\{\mathcal U\left[(g, \beta)\boldsymbol{\impl} (0_X,\alpha)\right]\otimes \mathcal G_{\beta}\Bigr\}\\
&=\mathcal U\left[(g, \beta)\boldsymbol{\impl} (0_X,\alpha)\right]\otimes \mathcal G_{\beta}.\\
\end{align*}
On the other hand, using the fact that  $\rho\leqslant \alpha$, we conclude that \linebreak $ \rho\vartriangleright \beta\leqslant \alpha \vartriangleright \beta$ which implies $(g,\beta)\boldsymbol{\impl} (0_X,\alpha)\preccurlyeq (g,\beta)\boldsymbol{\impl} (0_X,\rho)$.
\newline
Therefore,
\[
\mathcal U\left((g,\beta)\boldsymbol{\impl} (0_X,\alpha)\right)\leqslant\mathcal U\left((g,\beta)\boldsymbol{\impl} (0_X,\rho)\right),
\]
i.e.
\[
\mathcal G_{\beta}\leqslant\left[\mathcal U\left((g,\beta)\boldsymbol{\impl} (0_X,\rho)\right)\impl \bot\right].
\]
Now we invoke the residuation property of $(L,\leqslant,\otimes)$ to obtain
\[
\hat{\mathcal U}(0_X,\alpha)=\Bigl\{\mathcal U\left[(g, \beta)\boldsymbol{\impl} (0_X,\alpha)\right]\otimes \mathcal G_{\beta}\Bigr\}=\bot.
\]
Now we must show that $\hat{\mathcal U}$ is an $L$-fuzzy ultrafilter on $X$. In fact, since
\[
\hat{\mathcal U}(f,\alpha)=\mathcal U(f, \alpha)\bigvee\Bigl\{\mathcal U\left[(g, \beta)\boldsymbol{\impl} (f,\alpha)\right]\otimes \mathcal G_{\beta}\Bigr\},
\]
clearly $\mathcal U(f, \alpha)\leqslant\hat{\mathcal U}(f,\alpha),\,\,\,\ \forall (f,\alpha)\in L^X\times L$, but $\mathcal U$ is an $L$-fuzzy ultrafilter on $X$, therefore $\hat{\mathcal U}=\mathcal U$.
In this way
\begin{align*}
\mathcal U(g,\beta)&=\mathcal U(g, \beta)\bigvee\Bigl\{\mathcal U\left[(g, \beta)\boldsymbol{\impl} (g,\beta)\right]\otimes \mathcal G_{\beta}\Bigr\}\\
&=\mathcal U(g, \beta)\bigvee\Bigl\{\mathcal U(1_X, \bot)\otimes \mathcal G_{\beta}\Bigr\}\\
&=\mathcal U(g, \beta)\bigvee\Bigl\{\top\otimes \mathcal G_{\beta}\Bigr\}\\
&=\mathcal U(g, \beta)\lor \mathcal G_{\beta}.\\
\end{align*}
Therefore,
\[
\mathcal G_{\beta}=\left(\mathcal U\left[(g, \beta)\boldsymbol{\impl} (0_X,\rho)\right]\right)\impl\bot\leqslant \mathcal U(g, \beta),\,\,\,\ \forall (g,\beta)\in L^X\times L.
\]
From the last inequality and $(FF3')$ we obtain $(ii)$.\newline
$(ii)\Rightarrow(i)$\newline
We must verify that if
\[
\mathcal U(f,\alpha)=\left(\mathcal U\left[(f,\alpha)\boldsymbol{\impl}(0_X,\rho)\right]\right)\impl\bot,\,\ \text{for all}\,\ (f,\alpha)\in L^X\times L,\,\ \text {for all}\,\ \rho\leqslant\alpha\,\ \text{in} \,\ L
\]
then $\mathcal U$ is an $L$-fuzzy ultrafilter on $X$.\newline
Suppose  $\mathcal U\leqslant\hat{\mathcal U}$, then
\[
\left(\left(\hat{\mathcal U}\left[(f,\alpha)\boldsymbol{\impl}(0_X,\rho)\right]\right)\impl\bot\right)\leqslant \left(\left(\mathcal U\left[(f,\alpha)\boldsymbol{\impl}(0_X,\rho)\right]\right)\impl\bot\right),
\]
therefore $\hat{\mathcal U}\leqslant \mathcal U$, consequently $\mathcal U$ is an $L$-fuzzy ultrafilter on $X$.
\end{proof}
\begin{$^{*}$proposition}\label{di} 
Let $\phi: X\rightarrow Y$ be a map and let $\mathcal F:  L^X\times L\rightarrow L$ be  an $L$-fuzzy filter on $X$. Then
\begin{enumerate}
\item The map $\phi_{\mathcal F}^{\rightarrow}: L^Y\times L\rightarrow L$ defined by $\phi_{\mathcal F}^{\rightarrow}(g,\beta)= \mathcal F(g\circ\phi,\beta),\,\,\ \forall (g,\beta)\in L^Y\times L$ is an $L$-fuzzy filter on $Y$.
\item The map $\phi_{\mathcal U}^{\rightarrow}: L^Y\times L\rightarrow L$ is an $L$-fuzzy ultrafilter on $Y$, whenever $\mathcal U$ will be an $L$-fuzzy ultrafilter on $X$
\end{enumerate}
\end{$^{*}$proposition}
\begin{proof}
$(1).$ We must show $\phi_{\mathcal F}^{\rightarrow}$ satisfies the axioms of an $L$-fuzzy filter: In fact,
\begin{enumerate}
\item[(FF0)] $\phi_{\mathcal F}^{\rightarrow}(1_Y, \beta)= \mathcal F(1_Y\circ\phi,\beta)=\mathcal F(1_X,\beta)=\top,\,\,\ \forall \beta\in L$.
\item[(FF1)] $(f,\alpha)\preccurlyeq (g,\beta)\Rightarrow  \phi_{\mathcal F}^{\rightarrow}(f, \alpha)=\mathcal F(f\circ\phi,\alpha)
\leqslant \mathcal F(g\circ\phi, \beta)=\phi_{\mathcal F}^{\rightarrow}(g, \beta)$,
\newline since $(f\circ\phi,\alpha)\preccurlyeq (g\circ\phi,\beta) $.
\item[(FF2)] $\phi_{\mathcal F}^{\rightarrow}(f, \alpha)\otimes \phi_{\mathcal F}^{\rightarrow}(g, \beta)=\mathcal F(f\circ\phi,\alpha)\otimes \mathcal F(g\circ\phi,\beta)\leqslant \mathcal F((f\otimes g)\circ \phi, \alpha\lor \beta)= \phi_{\mathcal F}^{\rightarrow}(f\otimes g,\alpha\lor \beta)$,\newline
since $(f\circ\phi)\otimes (g\circ\phi)=[(f\otimes g)\circ \phi]$.
\item[(FF3)] $\phi_{\mathcal F}^{\rightarrow}(0_Y, \alpha)= \mathcal F(0_Y\circ\phi, \alpha)=\mathcal F(0_X, \alpha)=\bot,\,\,\,\ \forall \alpha\in L$.
\end{enumerate}
$(2).$ Let $\mathcal U:  L^X\times L\rightarrow L$ be  an $L$-fuzzy ultrafilter on $X$, let $(g,\beta)\in L^Y\times L$ and let $\alpha \in L$, then
\begin{align*}
\phi_{\mathcal U}^{\rightarrow}(g,\beta)&= \mathcal U(g\circ\phi,\beta)\\
&= \mathcal U\left[(g\circ\phi,\beta)\boldsymbol{\impl} (0_X,\alpha)\right]\impl \bot\\
&= \mathcal U[(g\circ\phi)\impl 0_X,\alpha\vartriangleright\beta]\impl \bot\\
&= \mathcal U[(g\circ\phi)\impl (0_Y\circ \phi),\alpha\vartriangleright\beta]\impl \bot\\
&= \mathcal U[(g\impl 0_Y)\circ \phi,\alpha\vartriangleright\beta]\impl \bot\\
&= \phi_{\mathcal U}^{\rightarrow}[g\impl 0_Y,\alpha\vartriangleright\beta]\impl \bot\\
&= \phi_{\mathcal U}^{\rightarrow}[(g,\beta)\boldsymbol{\impl} (0_Y,\alpha)]\impl \bot\\.
\end{align*}
We conclude from \ref{UF} that $\phi_{\mathcal U}^{\rightarrow}: L^Y\times L\rightarrow L$ is an $L$-fuzzy ultrafilter on $Y$.
\end{proof}
\begin{proposition}\label{ii}
Let $\phi: X\rightarrow Y$ be a surjective map and let\linebreak $\mathcal F:  L^Y\times L\rightarrow L$ be  an $L$-fuzzy filter on $Y$. Then the map $\phi_{\mathcal F}^{\leftarrow}: L^X\times L\rightarrow L$ defined by $\phi_{\mathcal F}^{\leftarrow}(f,\alpha)= \bigvee\left\{\mathcal F(g,\beta)\mid (g\circ\phi,\beta)\preccurlyeq (f, \alpha )  \right\},\,\,\ \forall (f,\alpha)\in L^X\times L$, is an $L$-fuzzy filter on $X$.
\end{proposition}
\begin{proof}
We must show $\phi_{\mathcal F}^{\leftarrow}$ satisfies the axioms of an $L$-fuzzy filter:
\begin{enumerate}
\item[(FF0)]
\begin{align*}
\phi_{\mathcal F}^{\leftarrow}(1_X, \alpha)&=\bigvee\left\{\mathcal F(g,\beta)\mid (g\circ\phi,\beta)\preccurlyeq (1_X, \alpha )  \right\}\\
&=\bigvee_{g\in L^Y}\bigwedge_{\beta \in L}\left\{\mathcal F(g,\beta)\mid g\circ \phi\leqslant 1_X,\,\  \alpha \leqslant \beta \right\}\\
&=\mathcal F(1_Y, \alpha)\\
&=\top,\,\,\ \forall \alpha \in L.
\end{align*}
\item[(FF1)] $(f,\alpha)\preccurlyeq (g,\beta)$ implies
\begin{align*}
\phi_{\mathcal F}^{\leftarrow}(f, \alpha)&=\bigvee\left\{\mathcal F(h,\delta)\mid (h\circ\phi,\delta)\preccurlyeq (f, \alpha )  \right\}\\
&=\bigvee_{h\in L^Y}\bigwedge_{\delta \in L}\left\{\mathcal F(h,\delta)\mid h\circ \phi\leqslant f,\,\  \alpha \leqslant \delta \right\}\\
&\leqslant \bigvee_{h\in L^Y}\bigwedge_{\delta \in L}\left\{\mathcal F(h,\delta)\mid h\circ \phi\leqslant g,\,\  \beta \leqslant \delta \right\}\\
&=\bigvee\left\{\mathcal F(h,\delta)\mid (h\circ\phi,\delta)\preccurlyeq (g, \beta ) \right\}\\
&=\phi_{\mathcal F}^{\leftarrow}(g\, \beta).
\end{align*}
\item[(FF2)]
\begin{align*}
&\phi_{\mathcal F}^{\leftarrow}(f, \alpha)\otimes \phi_{\mathcal F}^{\leftarrow}(g, \beta)=\\
&=\bigvee\left\{\mathcal F(h,\delta)\mid (h\circ\phi,\delta)\preccurlyeq (f, \alpha )  \right\}\otimes \bigvee\left\{\mathcal F(j,\eta)\mid (j\circ\phi,\eta)\preccurlyeq (g, \beta )  \right\}\\
&\leqslant \bigvee\ \{\mathcal F(h,\delta)\otimes F(j,\eta)\mid ((h\otimes j)\circ\phi,\delta\lor \eta)\preccurlyeq (f\otimes g, \alpha \lor\beta)  \}\\
&\leqslant \bigvee\left\{\mathcal F(h\otimes j,\delta\lor\eta)\mid ((h\otimes j)\circ\phi,\delta\lor \eta)\preccurlyeq (f\otimes g, \alpha \lor\beta)  \right\}\\
&=\phi_{\mathcal F}^{\leftarrow}(f\otimes g , \alpha\lor\beta).
\end{align*}
\item[(FF3)]
\begin{align*}\phi_{\mathcal F}^{\leftarrow}(0_X, \alpha)
&= \bigvee\left\{\mathcal F(h,\delta)\mid (h\circ\phi,\delta)\preccurlyeq (0_X, \alpha )  \right\}\\%
&=\bigvee_{h\in L^Y}\bigwedge_{\delta \in L}\left\{\mathcal F(h,\delta)\mid h\circ \phi\leqslant 0_X,\,\  \alpha \leqslant \delta \right\}\\
&=\mathcal F(0_Y, \alpha),\,\,\,\ \text{since $\phi$ is surjective} \\
&=\bot.
\end{align*}
\end{enumerate}
\end{proof}
\begin{definition}[\bf{$L$-fuzzy neighborhood system}]
Let $X$ be a set.  A map $\mathcal N: X \rightarrow L^{(L^X\times L)}$ is called an $L$-fuzzy neighborhood system on $X$ if and only if, for each $p\in X$, the mapping  $\mathcal N_{p}: L^X\times L \rightarrow L$ satisfies the following axioms:
\begin{enumerate}
\item[($N_0$)] $\mathcal N_p(1_X, \alpha) =\top$.
\item[($N_1$)] $(f,\alpha)\preccurlyeq (g,\beta)\Rightarrow  \mathcal N_p(f, \alpha)\leqslant \mathcal N_p(g, \beta)$.
\item[($N_2$)] $\mathcal N_p(f, \alpha)\otimes \mathcal N_p(g, \beta)\leqslant \mathcal N_p(f\otimes g, \alpha\lor \beta)$.
\item[($N_3$)] $\mathcal N_p(f, \alpha)\leqslant f(p)$.
\item[($N_4$)] $\mathcal N_p(f, \alpha)\leqslant\bigvee\left\{\mathcal N_p(g, \beta)\mid (f,\alpha)\preccurlyeq (g,\beta)\,\  \text{and}\,\ g(q)\leqslant\mathcal N_q(f, \alpha),\,\ \forall q\in X\right\}$.
\end{enumerate}
\end{definition}
\begin{definition}[\bf{$L$-fuzzy interior operator}]
Let $X$ be a set.  A map $\mathcal I:L^X\times L \rightarrow L^X$ is called an $L$-fuzzy interior operator on $X$ if and only if $\mathcal I$ satisfies the following conditions:
\begin{enumerate}
\item[($I_0$)] $\mathcal I(1_X, \alpha) =1_X,\,\,\ \forall \alpha \in L$.
\item[($I_1$)] $(f,\alpha)\preccurlyeq (g,\beta)\Rightarrow  \mathcal I(f, \alpha)\leqslant \mathcal I(g, \beta)$.
\item[($I_2$)] $\mathcal I(f, \alpha)\otimes \mathcal I(g, \beta)\leqslant \mathcal I(f\otimes g, \alpha\lor \beta)$.
\item[($I_3$)] $\mathcal I(f, \alpha)\leqslant f$.
\item[($I_4$)] $\mathcal I(f, \alpha)\leqslant \mathcal I\bigl(\mathcal I(f, \alpha)\bigr)$.
\item[($I_5$)] $\mathcal I(f, \bot)= f$.
\item[($I_6$)] If $\emptyset\ne K\subseteq L$,\ $\mathcal I(f, \alpha)= f^{0}\,\,\ \forall \alpha \in K$, then $I(f, \bigvee K)= f^{0}$.
\end{enumerate}
\end{definition}
\begin{lemma}(cf. \cite{HS})\label{int} 
Given an $L$-fuzzy topology $\mathcal T:L^X\rightarrow L$ on a set $X$, the mapping $\mathcal I_{\mathcal T}:L^X\times L\rightarrow L^X$ defined by
\[
\mathcal I(f,\alpha):=\mathcal I_{\mathcal T}(f,\alpha)=\bigvee \{u\in L^X\mid (u,\mathcal T(u))\preccurlyeq (f, \alpha)\},\,\,\,\ \forall (f,\alpha)\in L^X\times L
\]
is an $L$-fuzzy interior operator on $X$.
\end{lemma}
\begin{lemma}(cf. \cite{HS}) 
Every $L$-fuzzy interior operator  $\mathcal I:L^X\times L\rightarrow L^X$ induces, for each $p\in X$, an $L$-fuzzy  neighborhood system $\mathcal N^{\mathcal I}_{p}:L^X\times L\rightarrow L$ defined by
\[
\mathcal N_{p}(f,\alpha):=\mathcal N^{\mathcal I}_{p}(f,\alpha)=[\mathcal I(f,\alpha)](p).
\]
\end{lemma}
\begin{proposition}\label{cont}
Let $(X,\mathcal T)$ and $(Y,\sigma)$ be a pair of $L$-fuzzy topological spaces, let $\phi: X\rightarrow Y$ be an $L$-fuzzy continuous surjective  map, let\linebreak $\mathcal N_{p}: L^X\times L \rightarrow L$ be the $L$-fuzzy neighborhood system of a point $p\in X$ induced by $\mathcal T$, and let $\mathcal N_{\phi(p)}: L^Y\times L \rightarrow L$ be the corresponding  $L$-fuzzy neighborhood system of  $\phi(p)$ in $Y$ induced by $\sigma$. Then
\[
\mathcal N_{\phi(p)}\leqslant \phi^{\rightarrow}{(\mathcal N_{p})}
\]
where $\phi^{\rightarrow}{(\mathcal N_{p})}: L^Y\times L\rightarrow L$\  is defined by
\[
\phi^{\rightarrow}{(\mathcal N_{p})}(g,\beta)= \mathcal N_{p}(g\circ\phi,\beta),\,\,\ \forall (g,\beta)\in L^Y\times L.
\]
\end{proposition}
\begin{proof}
From lemma \ref{int} we have that
\[
\mathcal I_{\sigma}(g,\beta)=\bigvee \{u\in L^Y\mid (u,\sigma(u))\preccurlyeq (g, \beta)\},\,\,\,\ \forall (g,\beta)\in L^Y\times L
\]
is an $L$-fuzzy interior operator on $X$.
\noindent
The $L$-fuzzy continuity  of $\phi$ implies that
\[
\sigma(u)\leqslant \mathcal T(u\circ\phi),\,\,\,\,\,\,\ \forall u\in L^Y.
\]
We therefore have that
\[
\{u\in L^Y\mid u\leqslant g,\,\ \beta\leqslant \sigma(u)\}\subseteq \{v\in L^Y\mid v\leqslant g,\,\ \beta\leqslant \mathcal T(v\circ\phi)\},
\]
which implies
\[
\bigvee\{u\in L^Y\mid u\leqslant g,\,\ \beta\leqslant \sigma(u)\}\leqslant\bigvee \{v\in L^Y\mid v\leqslant g,\,\ \beta\leqslant \mathcal T(v\circ\phi)\}.
\]
In other words,
\[
\mathcal I_{\sigma}(g,\beta)\leqslant\bigvee \{v\in L^Y\mid v\leqslant g,\,\ \beta\leqslant \mathcal T(v\circ\phi)\}.
\]
If $\omega=\bigvee \{v\in L^Y\mid v\leqslant g,\,\ \beta\leqslant \mathcal T(v\circ\phi)\}$, then
\[
\mathcal N_{\phi(p)}(g,\beta)=\left[\mathcal I_{\sigma}(g,\beta)\right]\left(\phi(p)\right)\leqslant\omega\left(\phi(p)\right)=\left(\omega\circ\phi\right)(p),
\]
and so
\begin{equation}\label{1}
\mathcal N_{\phi(p)}(g,\beta)\leqslant\left(\omega\circ\phi\right)(p).
\end{equation}
On the other hand, it follows from $\omega\circ\phi\leqslant g\circ \phi$, and
\[
\mathcal N_{p}(g\circ \phi,\beta)=\left[ \mathcal I_{\mathcal T}(g\circ \phi,\beta)\right](p)
=\bigvee \{u\in L^X\mid u\leqslant g\circ \phi,\,\ \beta\leqslant \mathcal T(u)\}(p),
\]
and $\omega\circ\phi\in\{u\in L^X\mid u\leqslant g\circ \phi,\,\ \beta\leqslant \mathcal T(u)\}$, that
\begin{align*}
(\omega\circ\phi)(p)&\leqslant \bigvee\{u\in L^X\mid u\leqslant g\circ \phi,\,\ \beta\leqslant \mathcal T(u)\}(p)\\
&=\left[ \mathcal I_{\mathcal T}(g\circ \phi,\beta)\right](p)=\mathcal N_{p}(g\circ \phi,\beta)=\phi_{(\mathcal N_{p})}^{\rightarrow}(g,\beta),
\end{align*}
and so
\begin{equation}\label{2}
(\omega\circ\phi)(p)\leqslant \phi_{(\mathcal N_{p})}^{\rightarrow}(g,\beta).
\end{equation}
Finally, from \ref{1} and \ref{2}, we conclude that
\[
\mathcal N_{\phi(p)}\leqslant \phi^{\rightarrow}{(\mathcal N_{p})}.
\]
\end{proof}
\begin{definition}[\bf{Adherent point}] 
Let $(X,\tau)$ be an $L$-fuzzy topological space and $\mathcal N: X\rightarrow L^{(L^X\times L)}$ the corresponding  $L$-fuzzy neighborhood system. Further, let $\mathcal F:  L^X\times L\rightarrow L$ be an $L$-fuzzy filter on $X$. A point $p\in X$ is called an adherent point of $\mathcal F$ iff there exists a further  $L$-fuzzy filter $\mathcal G$  on $X$ provided with the following properties
\begin{enumerate}
\item[(i)] $\mathcal G(\bot\otimes\top).1_X, \alpha )\leqslant \bot\otimes\top,\,\,\,\ \forall \alpha \in L$.
\item[(ii)] $\mathcal N_p\leqslant \mathcal G$ and $\mathcal F\leqslant \mathcal G$.
\end{enumerate}
\end{definition}
\begin{definition}
An  $L$-fuzzy topological space $(X, \tau)$ is compact iff each  $L$-fuzzy filter on $X$ has at least  one adherent point.
\end{definition}
\begin{proposition}
Let $(X, \tau)$ be a compact $L$-fuzzy topological space, let $(Y, \sigma)$ be an $L$-fuzzy  topological space and let $\phi: X\rightarrow Y$ be a surjective, $L$-fuzzy continuous map. Then $(Y, \sigma)$ is compact.
\end{proposition}
\begin{proof}
Let $\mathcal F:L^Y\times L\rightarrow L$ be an $L$-fuzzy filter on $Y$, then by \ref{ii}, $\phi_{\mathcal F}^{\leftarrow}$ is an  $L$-fuzzy filter on $X$. Since $(X, \tau)$ is compact,  $\phi_{\mathcal F}^{\leftarrow}$ has an adherent point $p\in X$, i.e. there exists an $L$-fuzzy filter $\mathcal G$ on $X$ with $\mathcal N_p\leqslant \mathcal G$ and $\phi_{\mathcal F}^{\leftarrow}\leqslant \mathcal G$. Now we form, using \ref{di}, the corresponding $L$-fuzzy filter image and we obtain from the surjectivity    of $\phi$ the following relations
\[
\mathcal F=\phi^{\rightarrow}\left(\phi^{\leftarrow}(\mathcal F)\right)\leqslant \phi^{\rightarrow}(\mathcal G)\,\,\,\,\ \text{and}\,\,\,\,\ \phi^{\rightarrow}(\mathcal N_p)\leqslant \phi^{\rightarrow}(\mathcal G).
\]
On the other hand, from \ref{cont} follows
\[
\mathcal N_{\phi(p)}\leqslant \phi^{\rightarrow}{(\mathcal N_{p})},
\]
which implies that
\[
\mathcal F\leqslant \phi^{\rightarrow}(\mathcal G)\,\,\,\,\ \text{and}\,\,\,\,\ \mathcal N_{\phi(p)}\leqslant\phi^{\rightarrow}(\mathcal G),
\]
hence $\phi(p)$ is an adherent point of $\mathcal F$.\newline
This completes the proof of the proposition.
\end{proof}
In order to get a proof of our  version  of Tychonoff's theorem, we need the following lemmas:
\begin{lemma}\label{iii}
Let $\mathfrak F=\{(X_{\lambda},\mathcal T_{\lambda})\mid \lambda\in I\}$ be  a non-empty family of $L$-fuzzy topological spaces and let $(X,\mathcal T)$ be its $L$-fuzzy topological product. Then for each $ p\in X$, the mapping $\mathcal N_{p}: L^X\times L \rightarrow L$ defined by
\[
\mathcal N_{p}(f,\alpha)=\bigvee \left\{ \bigotimes_{\lambda\in I}\mathcal N_{p_{\lambda}}(h_{\lambda},\alpha)\mid h\in \Gamma_{f},\,\ \alpha\leqslant \bigotimes_{\lambda\in I}\tau_{\lambda}(h_{\lambda}) \right\},
\]
where
\[
\Gamma_f =\{ \mu \in \prod_{\lambda \in I}L^{X_{\lambda}}\mid \mu_{\lambda} =1_{X_{\lambda}} \text{\scriptsize{for all but finitely many indices $\lambda$, and}}\ \bigotimes_{\lambda \in I}(\mu_{\lambda}\circ\pi_{\lambda})\leqslant f\},
\]
is an $L$-fuzzy neighborhood system on $X$.
\end{lemma}
\begin{proof}
We must show that $\mathcal N_{p}$ satisfies the axioms of an $L$-fuzzy neighborhood system:
\newline
$(N_0).$ Let\ $1_{\Delta}\in\prod_{\lambda \in I}L^{X_{\lambda}}$\ defined by\ $(1_{\Delta})_{\lambda}=1_{X_{\alpha}}$,\ for each\ $\lambda \in I$.\ Then\ $1_{\Delta}\in \Gamma_{1_{X}}$. Since
\begin{enumerate}
\item $ \tau_{\lambda}(1_{X_{\lambda}})=\top$\ for all $\lambda\in I$, and so $\bigotimes_{\lambda\in I}\tau_{\lambda}(1_{X_{\lambda}})=\bigotimes_{\lambda\in I}\top=\top$,
\item $\mathcal N_{p_{\lambda}}(1_{X_{\lambda}},\alpha)=\top $\ for all $\lambda\in I$, and $\bigotimes_{\lambda\in I}\mathcal N_{p_{\lambda}}(1_{X_{\lambda}},\alpha)=\bigotimes_{\lambda\in I}\top=\top$,
then
\[
\mathcal N_{p}(1_X,\alpha)=\bigvee \left\{ \bigotimes_{\lambda\in I}\mathcal N_{p_{\lambda}}(h_{\lambda},\alpha)\mid h\in \Gamma_{1_X},\,\ \alpha\leqslant \bigotimes_{\lambda\in I}\tau_{\lambda}(h_{\lambda}) \right\}=\top.
\]
\end{enumerate}
We therefore have that $\mathcal N_{p}(1_X, \alpha)=\top. $
\newline
$(N_1).$ It is our purpose to show that $\mathcal N_{p}(f, \alpha)\leqslant\mathcal N_{p}(g, \beta)$, whenever $f\leqslant g$ in $L^X$, and $\beta \leqslant \alpha$ in $L$.
\newline
Clearly $f\leqslant g$ implies that $\Gamma_f\subseteq \Gamma_g$.
Consequently,
\begin{align*}
&\bigvee \left\{ \bigotimes_{\lambda\in I}\mathcal N_{p_{\lambda}}(h_{\lambda},\alpha)\mid h\in \Gamma_{f},\,\ \alpha\leqslant \bigotimes_{\lambda\in I}\tau_{\lambda}(h_{\lambda}) \right\}\\
\leqslant &\bigvee \left\{ \bigotimes_{\lambda\in I}\mathcal N_{p_{\lambda}}(h_{\lambda},\beta)\mid h\in \Gamma_{g},\,\ \beta\leqslant \bigotimes_{\lambda\in I}\tau_{\lambda}(h_{\lambda}) \right\}.
\end{align*}
This verifies that $\mathcal N_{p}(f, \alpha)\leqslant\mathcal N_{p}(g, \beta)$.
\newline
$(N_2).$ We must show that $\mathcal N_{p}(f, \alpha)\otimes\mathcal N_{p}(g, \beta)\leqslant \mathcal N_{p}(f\otimes g, \alpha\lor\beta)$.
\newline
Since\ $\mathcal N_{p}(f,\alpha)=\bigvee \left\{ \bigotimes_{\lambda\in I}\mathcal N_{p_{\lambda}}(h_{\lambda},\alpha)\mid h\in \Gamma_{f},\,\ \alpha\leqslant \bigotimes_{\lambda\in I}\tau_{\lambda}(h_{\lambda}) \right\}$\ and
\linebreak
$\mathcal N_{p}(g,\beta)=\bigvee \left\{ \bigotimes_{\lambda\in I}\mathcal N_{p_{\lambda}}(j_{\lambda},\beta)\mid j\in \Gamma_{g},\,\ \beta\leqslant \bigotimes_{\lambda\in I}\tau_{\lambda}(j_{\lambda}) \right\}$, and since $\otimes$ is distributive over non-empty joins, we have that
\begin{align*}
&\mathcal N_{p}(f, \alpha)\otimes\mathcal N_{p}(g, \beta)\\
&=\bigvee \left\{ \bigotimes_{\lambda\in I}\mathcal N_{p_{\lambda}}(h_{\lambda},\alpha)\otimes \bigotimes_{\lambda\in I}\mathcal N_{p_{\lambda}}(j_{\lambda},\beta)\mid h\in \Gamma_{f}, j\in\Gamma_{g},\ \alpha\lor \beta\leqslant \bigotimes_{\lambda\in I}\tau_{\lambda}(h_{\lambda}\otimes j_{\lambda})  \right\}\\
&=\bigvee \left\{ \bigotimes_{\lambda\in I}\left[\mathcal N_{p_{\lambda}}(h_{\lambda},\alpha)\otimes\mathcal N_{p_{\lambda}}(j_{\lambda},\beta)\right]\mid h\in \Gamma_{f}, j\in\Gamma_{g},\ \alpha\lor \beta\leqslant \bigotimes_{\lambda\in I}\tau_{\lambda}(h_{\lambda}\otimes j_{\lambda})  \right\}\\
&\leqslant \bigvee \left\{ \bigotimes_{\lambda\in I}\mathcal N_{p_{\lambda}}(h_{\lambda}\otimes j_{\lambda},\alpha\lor\beta)\mid h\otimes j\in \Gamma_{f\otimes g},\ \alpha\lor \beta\leqslant \bigotimes_{\lambda\in I}\tau_{\lambda}(h_{\lambda}\otimes j_{\lambda})  \right\}\\
&=\mathcal N_{p}(f\otimes g, \alpha\lor\beta).
\end{align*}
$(N_3).$ We must verify that $  N_{p}(f, \alpha)\leqslant f(p)$, for each $(f,\alpha)\in L^X\times L$.
\newline
Since $\mathcal N_{p_{\lambda}}$ is an $L$-fuzzy neighborhood system on $X_{\lambda}$, we have that\linebreak $\mathcal N_{p_{\lambda}}(h_{\lambda},\alpha)\leqslant h_{\lambda}(p_{\lambda})$, for all $\lambda\in I$, and for  all $(h_{\lambda},\alpha)\in L^{X_{\lambda}}\times L$. Now, let $h\in \Gamma_f$, then $h \in \prod_{\lambda \in I}L^{X_{\lambda}},\,\ h_{\lambda} =1_{X_{\lambda}}$\  for all but finitely many indices $\lambda\in I $, and $ \bigotimes_{\lambda \in I}(h_{\lambda}\circ\pi_{\lambda})\leqslant f$. Hence
\[
\bigotimes_{\lambda \in I}(h_{\lambda}\circ\pi_{\lambda})(p)=  \bigotimes_{\lambda \in I}(h_{\lambda})(p_{\lambda})\leqslant f(p).
\]
We therefore have that
\[
\bigotimes_{\lambda\in I}\mathcal N_{p_{\lambda}}(h_{\lambda},\alpha)\leqslant\bigotimes_{\lambda \in I}(h_{\lambda})(p_{\lambda})\leqslant f(p),
\]
which implies
\[
\mathcal N_{p}(f,\alpha)=\bigvee \left\{ \bigotimes_{\lambda\in I}\mathcal N_{p_{\lambda}}(h_{\lambda},\alpha)\mid h\in \Gamma_{f},\,\ \alpha\leqslant \bigotimes_{\lambda\in I}\tau_{\lambda}(h_{\lambda}) \right\}\leqslant f(p).
\]
$(N_4).$ To show that
\[
\mathcal N_p(f, \alpha)\leqslant\bigvee\left\{\mathcal N_p(g, \beta)\mid (f,\alpha)\preccurlyeq (g,\beta),\  g(q)\leqslant\mathcal N_q(f, \alpha),\ \forall q\in X\right\},
\]
we must verify the following:
\[
\mathcal N_{p_{\lambda}}(\omega_{\lambda}, \alpha)\leqslant\bigvee\left\{\mathcal N_{p_{\lambda}}(v_{\lambda}, \beta)\mid (\omega_{\lambda},\alpha)\preccurlyeq (v_{\lambda},\beta),\ v_{\lambda}(q_{\lambda})\leqslant\mathcal N_{q_{\lambda}}(\omega_{\lambda}, \alpha),\ \forall q_{\lambda}\in X_{\lambda}\right\}
\]
implies
\[
\mathcal N_p(f, \alpha)\leqslant\bigvee\left\{\mathcal N_p(g, \beta)\mid (f,\alpha)\preccurlyeq (g,\beta),\   g(q)\leqslant\mathcal N_q(f, \alpha),\ \forall q\in X\right\}.
\]
We observe from the hypothesis that, for $\omega\in \Gamma_f$\ and $v\in \Gamma_g$,
\begin{align*}
&\bigotimes_{\lambda\in I}\mathcal N_{p_{\lambda}}(\omega_{\lambda}, \alpha)\\
&\leqslant\bigvee\left\{\bigotimes_{\lambda\in I}\mathcal N_{p_{\lambda}}(v_{\lambda}, \beta)\mid (\omega_{\lambda},\alpha)\preccurlyeq (v_{\lambda},\beta),\ v_{\lambda}(q_{\lambda})\leqslant\mathcal N_{q_{\lambda}}(\omega_{\lambda}, \alpha),\ \forall q_{\lambda}\in X_{\lambda}\right\}\\
&\leqslant\bigvee\left\{\bigotimes_{\lambda\in I}\mathcal N_{p_{\lambda}}(v_{\lambda}\circ\pi_{\lambda}, \beta)\mid (\omega_{\lambda},\alpha)\preccurlyeq (v_{\lambda},\beta),\ v_{\lambda}(q_{\lambda})\leqslant\mathcal N_{q_{\lambda}}(\omega_{\lambda}, \alpha),\ \forall q_{\lambda}\in X_{\lambda}\right\}.\\
\end{align*}
We therefore have that
\[
\mathcal N_p(f, \alpha)\leqslant\bigvee\left\{\mathcal N_p(g, \beta)\mid (f,\alpha)\preccurlyeq (g,\beta),\   g(q)\leqslant\mathcal N_q(f, \alpha),\ \forall q\in X\right\}.
\]
\end{proof}
\begin{$^{*}$lemma}
Let $\mathfrak F=\{(X_{\lambda},\mathcal T_{\lambda})\mid {\lambda}\in I\}$ be  a non-empty family of compact $L$-fuzzy topological spaces and let $(X,\mathcal T)$ be its $L$-fuzzy topological product. Let $\mathcal U:L^X\times L\rightarrow L$ be an $L$-fuzzy ultrafilter on $X$, then the following statements are equivalent:
\begin{enumerate}
\item  $\mathcal U$ converges to an element $p=(p_{\lambda})_{\lambda\in I}$\,\ in $X$.
\item For each $\lambda\in I$,\,\ $\pi_{\lambda}^{\rightarrow}(\mathcal U)$ converges to $p_{\lambda}$, where $\pi_{\lambda};X\rightarrow X_{\lambda}$ is the $\lambda$th projection.
\end{enumerate}
\end{$^{*}$lemma}
\begin{proof}
$(1)\Rightarrow (2)$ is trivial.\newline
$(2)\Rightarrow (1).$ Suppose that $\pi_{\lambda}^{\rightarrow}(\mathcal U)$ converges to $p_{\lambda}$, for each $\lambda\in I$. Then $\mathcal N_{p_{\lambda}}\leqslant \pi_{\lambda}^{\rightarrow}(\mathcal U)$. Therefore
\begin{align*}
\bigotimes_{\lambda\in I}\mathcal N_{p_{\lambda}}(h_{\lambda},\alpha)&\leqslant \bigotimes_{\lambda\in I}\pi_{\lambda}^{\rightarrow}(\mathcal U)(h_{\lambda},\alpha)\\
&=\bigotimes_{\lambda\in I}\mathcal U(h_{\lambda}\circ \pi_{\lambda},\alpha)\\
&\leqslant\mathcal U\left(\bigotimes_{\lambda\in I}(h_{\lambda}\circ \pi_{\lambda}),\alpha \right)\\
&\leqslant\mathcal U\left(f,\alpha \right).\\
\end{align*}
It follows that
\[
\mathcal N_{p}(f,\alpha)=\bigvee \left\{ \bigotimes_{\lambda\in I}\mathcal N_{p_{\lambda}}(h_{\lambda},\alpha)\mid h\in \Gamma_{f},\,\ \alpha\leqslant \bigotimes_{\lambda\in I}\tau_{\lambda}(h_{\lambda}) \right\}\leqslant\mathcal U\left(f,\alpha \right).
\]
This establishes that $\mathcal U$ converges to an element $p=(p_{\lambda})_{\lambda\in I}$\,\ in $X$.
\end{proof}
As a consequence of the previous lemmas, we get the main result:
\begin{$^{*}$theorem}
Let $\mathfrak F=\{(X_{\lambda},\mathcal T_{\lambda})\mid \lambda\in I\}$ be  a non-empty family of $L$-fuzzy topological spaces. Then the following assertions are equivalent:
\begin{enumerate}
\item $(X_{\lambda},\mathcal T_{\lambda})$ is compact, for all  $\lambda\in I$.
\item The  $L$-fuzzy topological product $(X_0, \mathcal T_0)$ of $\mathfrak F$ in the  category $\mathbf {L\mbox{-}FTop}$ is compact.
\end{enumerate}
\end{$^{*}$theorem}


\begin{thebibliography}{99}

\bibitem{AHS} J.Ad\'amek, H. Herrlich, G. E. Strecker, {\it Abstract And Concrete Categories: The Joy Of Cats}, Wiley Interscince Pure And Applied Mathematics, John Wiley \& Sons (Brisbane/Chisester/New York/Singapure/Toronto), 1990.

\bibitem{GB} G. Birkhoff, {\it Lattice Theory}, American Mathematical Society Colloquium Publications, Volume \textbf{XXV},  The American Mathematical Society (Providence, Rhode Island), 1940.

\bibitem{NB} N. Bourbaki, {\it General Topology, Part 1}, Addison-Wesley Publishing Company, Advanced Book Program Reading, Massachusetts, 1966.

\bibitem{Ch58} C.C. Chang, {\it Algebraic analysis of many valued logics}, Trans. Amer. Math. Soc. \textbf{88}(1958), 467--490.

\bibitem{Ch59} \bysame, {\it A new proof of the completeness of the {\l}ukasiewicz axiom},  Trans. Amer. Math. Soc. \textbf{93}(1959), 74--80.

\bibitem{Ho91} U. H{\"o}hle, {\it Monoidal closed categories, weak topoi and generalized  logics}, Fuzzy Sets and Systems \textbf{42}(1991), 15--35.

\bibitem{Ho92} \bysame, {\it ${M}$-valued sets and sheaves over integral commutative ${\rm cl}$-monoids}, Applications of category theory to fuzzy subsets (Linz, 1989),  Kluwer Acad. Publ., Dordrecht, 1992, 33--72.

\bibitem{Ho94} \bysame, {\it Commutative, residuated $l$-monoids}, Non-classical logics and  their applications to fuzzy subsets (Linz, 1992), Kluwer Acad. Publ.,  Dordrecht, 1995, 53--106.

\bibitem{HS} U. H\"ohle, A. \v{S}ostak, {\it Axiomatic fundations of fixed-basis fuzzy topology}, Chapter 3 in: U. H\"ohle, S. E. Rodabaugh, eds, Mathematics Of Fuzzy Sets: Logic, Topology And Measure Theory, The Handbook of Fuzzy Sets Series, Volume \textbf{3}(1999), Kluwer Academic Publisher (Boston/Dordrecht/London), 123--272.

\bibitem{Jhst} P. T. Johnstone, {\it Stone spaces}, Cambridge University Press (Cambridge), 1982.

\bibitem{SER} S. E. Rodabaugh, {\it Powerset Operator Foundations For Poslat Fuzzy Set Theories and Topologies}, Chapter 2 in: U. H\"ohle, S. E. Rodabaugh, eds, Mathematics Of Fuzzy Sets: Logic, Topology And Measure Theory, The Handbook of Fuzzy Sets Series, Volume \textbf{3}(1999), Kluwer Academic Publisher (Boston/Dordrecht/London), 91--116.

\bibitem{AS} A P. \v{S}ostak, {\it Fuzzy functions and an extension of the category\ $L$-Top of Chang-Goguen\ $L$-topological spaces}, Proceedings of the Ninth Prague Topological Symposium, Contributed papers from the symposium held in Prague, Czech Republic, August 19--25, 2001, 271--294, Topology Atlas, Toronto, 2002; {\tt arXiv:math.GN/0204139}.

\bibitem{SW} S. Willard, {\it General Topology}, Addison-Wesley Publishing Company (Reading, Massachusetts/Menlo Park, California/London/Amsterdam/Don Mills, Ontario/Sydney), 1970.

\end{thebibliography}
\end{document}